\documentclass[conference,letter]{ieeeconf}
\usepackage{cite}
\usepackage{amsmath,amssymb,amsfonts}
\usepackage{algorithmic}
\usepackage{graphicx}
\usepackage{textcomp}
\usepackage{xcolor}

\usepackage{tikz}
\usetikzlibrary{backgrounds,matrix,fit}

\usepackage{epstopdf}
\usepackage[caption=false]{subfig}

\newcommand*\circled[1]{\tikz[baseline=(char.base)]{
\node[shape=circle, draw, inner sep=0.6pt] (char) {#1};}
}
\newcommand\kronF[2]{#1^{\circled{\tiny{#2}}}}

\graphicspath{./figures}

\begin{document}

\title{Sum of squares approximations to energy functions}

\author{\authorblockN{Hamza Adjerid}
\authorblockA{\textit{Department of Mathematics} \\
\textit{Virginia Tech}\\
Blacksburg, VA USA \\
ORCID 0009-0000-9025-0764}
\and
\authorblockN{Jeff Borggaard}
\authorblockA{\textit{Department of Mathematics} \\
\textit{Virginia Tech}\\
Blacksburg, VA USA \\
ORCID 0000-0002-4023-7841}
}

\maketitle

\begin{abstract}
Energy functions offer natural extensions of controllability and observability Gramians to nonlinear systems, enabling various applications such as computing reachable sets, optimizing actuator and sensor placement, performing balanced truncation, and designing feedback controllers. However, these extensions to nonlinear systems depend on solving Hamilton-Jacobi-Bellman (HJB) partial differential equations, which are infeasible for large-scale systems. Polynomial approximations are a viable alternative for modest-sized systems, but conventional polynomial approximations may yield negative values of the energy away from the origin. To address this issue, we explore polynomial approximations expressed as a sum of squares to ensure non-negative approximations. In this study, we focus on a reduced sum of squares polynomial where the coefficients are found through least-squares collocation---minimizing the HJB residual at sample points within a desired neighborhood of the origin. We validate the accuracy of these approximations through a case study with a closed-form solution and assess their effectiveness for controlling a ring of van der Pol oscillators with a Laplacian-like coupling term and discretized Burgers equation with source terms.
\end{abstract}


\section{Introduction}
Energy functions are pervasive in control theory.  They were developed to extend  controllability and observability Gramians to nonlinear systems for balancing~\cite{scherpen1993BalancingNonlinearSystems} and have since been generalized to the $H_\infty$ framework as the basis for nonlinear balanced truncation~\cite{scherpen1996Hinf}.  Since they can also serve as Lyapunov functions for closed-loop systems, other applications of energy functions include estimation of reachable sets and safe control algorithms~\cite{bansal2017HamiltonJacobiReachabilityBrief,wei2019SafeControlAlgorithms}, sensor and actuator placement~\cite{georges1995UseObservabilityControllability}, and control design cf.~\cite{haddad2018EnergybasedFeedbackControl}.  The computational challenge to the use of energy functions are their reliance on solutions to the Hamilton-Jacobi-Bellman (HJB) equations.  To avoid the curse-of-dimensionality, polynomial approximations have been introduced, e.g.~\cite{fujimoto2008ComputationNonlinearBalanced,krener2008ReducedOrderModeling,kramer2023nonlinear1}.  However, these approximations often exhibit negative energy away from the origin.

Sum of squares (SOS) approximations in control theory also has a rich history that arose as a result of Parillo~\cite{parillo2000structured} recasting many sum of squares problems as convex optimization problems and the development of efficient algorithms to solve them~\cite{boyd2004ConvexOptimization}. Their primary role in systems theory is to develop Lyapunov functions, which are guaranteed to be positive definite with negative definite derivatives along solutions and are useful for control synthesis~\cite{prajnaNonlinearControlSynthesis2004,papachristodoulouTutorialSumSquares2005}. 
It is natural to consider SOS approximations for energy functions, and this has been the subject of a number of studies, cf.~\cite{jarvis-wloszek2005ControlsApplicationsSum} and have specifically been used for energy functions for a policy iteration control algorithm in recent work~\cite{pakkhesal2022SumSquaresBased}.

The remainder of this paper is organized as follows.  In the next section, we provide a background on energy functions, a recent Kronecker product-based polynomial approximation that leads to efficient calculation of coefficients for modest systems, and a review of a typical sum of squares representation.  Then we introduce three SOS formulations.  The first embeds the Kronecker product polynomials within a higher degree polynomial to guarantee SOS.  A second couples the Kronecker product polynomial with a nonlinear least-squares collocation procedure to determine the higher degree terms in a SOS.  Finally, the third involves a modification of Parrilo's formulation and nonlinear least-squares collocation.  To aid the optimization, we
introduce a collocation strategy where we minimize residuals that are sampled over wider parameter sets.  This allows us to warm-start the optimization to compute the SOS approximation in the region of interest.  Numerical examples include a scalar test case where the analytic solution is available, a ring of van der Pol oscillators that are coupled through a Laplacian-like term, and an example generated by finite element discretization of Burgers equation.  These show the better approximations to the energy functions away from the origin when using the modified Parrilo formulation and the warm-start strategy.

\section{Background}

\subsection{Energy functions}
Consider the nonlinear, control affine input-output system
\begin{equation}
\label{eq:nl_system}
  \dot{{\bf x}}(t)={\bf f}({\bf x}(t))+{\bf Bu}(t), \qquad {\bf y}(t)={\bf Cx}(t)
\end{equation}
where ${\bf f}:\mathbb{R}^n\rightarrow \mathbb{R}^n$ is a smooth function, ${\bf B}\in \mathbb{R}^{n\times m}$,
and ${\bf C}\in \mathbb{R}^{p\times n}$.  The space of admissible controls are 
\begin{displaymath}
{\cal U}^- = L_2(-\infty,0;\mathbb{R}^m) \qquad \mbox{and} \qquad {\cal U}^+ = L_2(0,\infty;\mathbb{R}^m).
\end{displaymath}
For this class of systems, energy functions are a natural extension of $H_\infty$ control theory~\cite{vanderschaft1992L2gainAnalysisNonlinear}.  For any $\eta\leq 1$, we define
\begin{equation}
\label{eq:pastEnergy}
{\mathcal{E}}_\eta^{-}({\mathbf{x}}_0)  :=\min_{\substack{{\bf u} \in {\cal U}^- \\ {\bf x}(-\infty) = {\bf 0} \\ {\bf x}(0) = {\bf x}_0}} \ \frac{1}{2} \int_{-\infty}^{0} \eta\Vert {\mathbf{y}}(t) \Vert^2  +  \Vert {\mathbf{u}}(t) \Vert^2 {\rm{d}}t
\end{equation}
as the {\em past energy}.  The {\em future energy} is defined for three cases.  For $0<\eta\leq 1$, 
\begin{equation}
\label{eq:futureEnergy1}
{\mathcal{E}}_\eta^{+}({\mathbf{x}}_0)  :=\min_{\substack{{\bf u} \in {\cal U}^+ \\ {\mathbf{x}}(0) = {\mathbf{x}}_0 \\ {\mathbf{x}}(\infty) = {\bf 0}}} \ \frac{1}{2} \int_{0}^{\infty} \Vert {\bf y}(t) \Vert^2  +  
\frac{1}{\eta}\Vert {\bf u}(t) \Vert^2 {\rm{d}}t;
\end{equation}
for $\eta=0$,
\begin{equation}
\label{eq:controllability}
{\mathcal{E}}_\eta^{+}({\mathbf{x}}_0) := \frac{1}{2} \int \Vert {\bf y}(t) \Vert^2{\rm{d}}t, \quad {\mathbf{x}}(0)={\mathbf{x}}_0, \ \mbox{and}\ {\bf u}(t)\equiv {\bf 0};
\end{equation}
and for $\eta<0$,
\begin{equation}
\label{eq:futureEnergy3}
{\mathcal{E}}_\eta^{+}({\mathbf{x}}_0)  :=\max_{\substack{{\bf u} \in {\cal U}^+ \\ {\mathbf{x}}(0) = {\mathbf{x}}_0, \\  {\mathbf{x}}(\infty) = {\bf 0}}} \ \frac{1}{2} \int\displaylimits_{0}^{\infty} \Vert {\bf y}(t) \Vert^2  +  \frac{1}{\eta}\Vert {\bf u}(t) \Vert^2 {\rm{d}}t.
\end{equation}
These are often defined using an $L_2$-gain parameter $\gamma_0$ when $\eta = -\gamma_0^{-2}$ \cite{vanderschaft1992L2gainAnalysisNonlinear} or using an $H_\infty$-gain parameter $\gamma>0$ when $\eta = 1-\gamma^{-2}$ \cite{scherpen1996Hinf}.
In the $H_\infty$ setting, when we take the limit $\gamma\rightarrow 1$, $\eta=0$, we recover that ${\mathcal{E}}_0^{-}$ is the controllability energy and ${\mathcal{E}}_0^{+}$ is the observability energy.  Additionally, for the limit $\gamma\rightarrow \infty$, $\eta=1$, we recover the energy functions used in HJB balancing \cite{scherpen1994normalized}.  Finally, when $\eta=-1$, the forward and past energy functions agree.  These connections are detailed in \cite{kramer2023nonlinear1}.  

These energy functions can be characterized by HJB equations, cf.~\cite{scherpen1996Hinf}.  If $\bar{\mathcal E}$ is a solution to
\begin{equation} \label{eq:HJB-NLHinfty2}
 0  =  \frac{\partial \bar{\mathcal{E}}({\mathbf{x}})}{\partial {\mathbf{x}}} {\bf f}({\mathbf{x}}) + \frac{1}{2}  \frac{\partial \bar{\mathcal{E}}({\mathbf{x}})}{\partial {\mathbf{x}}} {\bf B} {\bf B}^\top   \frac{\partial^\top   \bar{\mathcal{E}}({\mathbf{x}})}{\partial {\mathbf{x}}}
 - \frac{\eta}{2} {\bf x}^\top   {\bf C}^\top    {\bf Cx}
\end{equation}
with $\bar{\mathcal{E}}({\bf 0}) = 0$ and ${\bf 0}$ is an asymptotically stable fixed point of 
\begin{equation}
  \dot{\bf x} =  - \left ( {\bf f}({\mathbf{x}}) +{\bf B}{\bf B}^\top   \frac{\partial^\top   \bar{\mathcal{E}}({\mathbf{x}})}{\partial {\mathbf{x}}} \right ), 
\end{equation}
then $\bar{\mathcal{E}}({\bf x})$ is the past energy function ${\mathcal{E}}_\eta^{-}({\mathbf{x}})$ from \eqref{eq:pastEnergy}.

Likewise, if $\tilde{\mathcal{E}}$ is a solution to
\begin{equation} \label{eq:HJB-NLHinfty1}
0  = \frac{\partial \tilde{\mathcal{E}}({\mathbf{x}})}{\partial {\mathbf{x}}} {\bf f}({\mathbf{x}})
 - \frac{\eta}{2}  \frac{\partial \tilde{\mathcal{E}}({\mathbf{x}})}{\partial {\mathbf{x}}} {\bf B} {\bf B}^\top   \frac{\partial^\top   \tilde{\mathcal{E}}({\mathbf{x}})}{\partial {\mathbf{x}}}
 + \frac{1}{2}{\bf x}{\bf C}^\top   {\bf C}{\bf x}
\end{equation}
with $\tilde{\mathcal{E}}({\bf 0}) = 0$ and  ${\bf 0}$ is an asymptotically stable fixed point of
\begin{equation}
   \dot{\bf x} =  {\bf f}({\mathbf{x}}) - \eta{\bf B} {\bf B}^\top   \frac{\partial^\top   \tilde{\mathcal{E}}({\mathbf{x}})}{\partial {\mathbf{x}}}, 
\end{equation}
then this solution $\tilde{\mathcal E}({\bf x})$ is the future energy function ${\mathcal{E}}_\eta^{+}({\mathbf{x}})$.

\subsection{Polynomial approximations\label{sec:polynomial}}
Mesh-based approximations to Hamilton-Jacobi equations are prohibitive for systems with modest dimension.  Therefore approximations using global basis functions, typically polynomials, were soon introduced~\cite{albrekht1961optimal,lukes1969optimal,navasca2000solution,breiten2019taylor}.  Polynomial approximations to energy functions in the forms defined above (\ref{eq:pastEnergy})--(\ref{eq:futureEnergy3}) were developed for nonlinear balancing applications as well~\cite{fujimoto2008ComputationNonlinearBalanced}.  
In \cite{kramer2023nonlinear1}, the authors developed a scalable approach to approximate energy functions by polynomials when they are written using Kronecker products:  
\begin{align}
\label{eq:pastEpoly}
    \mathcal{E}_\eta^{-}(\mathbf{x})\approx&\cfrac{1}{2}\left( \mathbf{v}_2^\top \kronF{\mathbf{x}}{2}+\mathbf{v}_3^\top \kronF{\mathbf{x}}{3}+\cdots +\mathbf{v}_d^\top \kronF{\mathbf{x}}{d}\right)\\
\label{eq:futureEpoly}
\mathcal{E}_\eta^{+}(\mathbf{\mathbf{x}})\approx&\cfrac{1}{2}\left( \mathbf{w}_2^\top \kronF{\mathbf{x}}{2}+\mathbf{w}_3^\top \kronF{\mathbf{x}}{3}+\cdots +\mathbf{w}_d^\top \kronF{\mathbf{x}}{d}\right),
\end{align}
where $\kronF{\mathbf{x}}{d}$ is a Kronecker product involving $d$ copies of the vector $\mathbf{x}$, e.g. $\kronF{\mathbf{x}}{3} = \mathbf{x}\otimes\mathbf{x}\otimes\mathbf{x}$.  
The coefficients ${\bf v}_k$ and ${\bf w}_k$ involve structured linear systems that enable efficient computation.  For ${\bf M}\in\mathbb{R}^{q\times n}$ we define the \textit{k-way Lyapunov matrix} or a special \textit{Kronecker sum} \cite{benzi2017approximation} matrix as
\begin{displaymath}
    \mathcal{L}_k({\bf M})\! :=\! \underbrace{{\bf M} \otimes \ldots \otimes {\bf I}_n}_{k \  \text{times}} + \cdots + \underbrace{{\bf I}_n \otimes \ldots \otimes {\bf M}}_{k \ \text{times}} \in \mathbb{R}^{n^{k-1}q \times n^k},
\end{displaymath}
%
where ${\bf I}_n$ is the $n$-dimensional identity matrix. 
We present the results from \cite{kramer2023nonlinear1} for calculating the coefficients ${\bf v}_k$ and ${\bf w}_k$ when system (\ref{eq:nl_system}) is specialized to $\mathbf{f}({\bf x})\equiv {\bf Ax}+{\bf F}\kronF{\bf x}{2}$, and the pairs $({\bf A},{\bf B})$ and $({\bf A},{\bf C})$ are stabilizable and detectable. \newline

\textbf{Theorem 1:} (\!\!\cite[Th.~7]{kramer2023nonlinear1}) Consider system (\ref{eq:nl_system}) with the conditions stated above.  Given $\eta\leq 1$, and past energy function $\mathcal{E}^-_\eta(\mathbf{x})$ expanded with the
coefficients $\mathbf{v}_i,$  $i = 2, \ldots, d$  in (\ref{eq:pastEpoly}). Then, $\mathbf{v}_2 = {\tt vec} (\mathbf{V}_2)$,
where $\mathbf{V}_2$ is the symmetric positive definite solution to
the $\mathcal{H_\infty}$ Riccati equation
\begin{equation}\label{eq:12}
 0=\mathbf{A}^\top   \mathbf{V}_2+\mathbf{V}_2\mathbf{A}-\eta \mathbf{C}^\top   \mathbf{C}+\mathbf{V}_2\mathbf{B}\mathbf{B}^\top   \mathbf{V}_2.   
\end{equation}

Moreover, the coefficient vectors $\mathbf{v}_k = {\tt vec}(\mathbf{V}_k) \in \mathbb{R}^{n^k}$ for $2\leq k\leq d$ solve the linear systems
\begin{IEEEeqnarray*}{lCr}
    \mathcal{L}_k((\mathbf{A} + \mathbf{B}\mathbf{B}^\top\mathbf{V}_2  )^\top)\mathbf{v}_k= \\
    \ -\mathcal{L}_{k-1}(\mathbf{F}^\top  )\mathbf{v}_{k-1}-\cfrac{\eta}{4}\sum_{\substack{i,j>2\\ i+j=k+2}} ij\,{\tt vec}(\mathbf{V}_i^\top \mathbf{B}\mathbf{B}^\top \mathbf{V}_j).
\end{IEEEeqnarray*}

\textbf{Theorem 2:} (\!\!\cite[Th.~6]{kramer2023nonlinear1}) Consider system (\ref{eq:nl_system}) with the conditions stated above. Given $\eta\leq 1$ and future energy function $\mathcal{E}^+_\eta(\mathbf{x})$ expanded with the
coefficients ${\bf w}_i,$ $i = 2, \ldots, d$ in (\ref{eq:futureEpoly}). Then, $\mathbf{w}_2 = {\tt vec} (\mathbf{W}_2)$,
where $\mathbf{W}_2$ is the symmetric positive definite solution to
the $\mathcal{H_\infty}$ Riccati equation
\begin{equation}\label{eq:13}
 0=\mathbf{A}^\top  \mathbf{W}_2+\mathbf{W}_2\mathbf{A}+\mathbf{C}^\top  \mathbf{C}-\eta \mathbf{W}_2\mathbf{B}\mathbf{B}^\top \mathbf{W}_2.   
\end{equation}

Moreover, the coefficient vectors $\mathbf{w}_k = {\tt vec}(\mathbf{W}_k) \in \mathbb{R}^{n^k}$ for $2\leq k\leq d$ solve the linear systems 
\begin{IEEEeqnarray*}{lCr}
\mathcal{L}_k((\mathbf{A} -\eta \mathbf{B}\mathbf{B}^\top\mathbf{W}_2 )^\top)\mathbf{w}_k=\\
\ -\mathcal{L}_{k-1}(\mathbf{F}^\top )\mathbf{w}_{k-1}+\cfrac{\eta}{4}\sum_{\substack{i,j>2\\ i+j=k+2}} ij \, {\tt vec}(\mathbf{W}_i^\top  \mathbf{B}\mathbf{B}^\top  \mathbf{W}_j).  
\end{IEEEeqnarray*}

The proofs for both theorems can be found in \cite{kramer2023nonlinear1}. While solving the systems grows exponentially in $k$ and as a polynomial in $n$, the authors developed an efficient implementation taking advantage of the Kronecker structure of the systems based on the work in \cite{borggaard2019QQR,borggaard2021PQR,chen2019RecursiveBlockedAlgorithms}. 

As shown in \cite{kramer2023nonlinear1}, polynomial approximations are accurate near the origin. The main issue with these approximations is negativity away from the origin.  However, energy functions must be positive definite by definition. One way to overcome this issue is to propose function approximations that impose non-negative definiteness.

\subsection{Sum of squares}
A scalar valued function $f_{sos}(\mathbf{x})$ is a sum of squares if it can be written as 
\begin{displaymath}
  f_{\rm sos}(\mathbf{x})=\sum_{i=1}^N f_i^2(\mathbf{x}),
\end{displaymath}
where $\{f_i(\mathbf{x})\}$ can be any collection of generic functions and $N$ any integer. Polynomial sum of squares are much easier to work with
\begin{displaymath}
  f_{\rm sos}(\mathbf{x})=\sum_{i=1}^N p_i^2(\mathbf{x}),
\end{displaymath}
where each $p_i(\mathbf{x})$ is a polynomial. 
\newline

\textbf{Proposition}: (\!\!\cite{parillo2000structured}) A polynomial $p(\mathbf{x})$ of degree $2d$ is a SOS if and only if there exists a positive semi-definite matrix $\mathbf{Q}$ and vector  $\mathbf{z}(\mathbf{x})$ that contains monomials in $\mathbf{x}$ of degree $\leq$ d such that: \begin{equation}
\label{eq:psos}
p(\mathbf{x})=\mathbf{z}(\mathbf{x})^\top \mathbf{Q}\mathbf{z}(\mathbf{x}).
\end{equation}

For example~\cite{parillo2000structured}, by defining
$\mathbf{z}(\mathbf{x})=[x_1^2 \quad x_1x_2 \quad x_2^2]^\top$, we can rewrite
\begin{displaymath}
  p({\bf x})=2x_1^4+2x_1^3x_2-x_1^2x_2^2+5x_2^4
\end{displaymath}
in the form (\ref{eq:psos}) with the positive semi-definite matrix
\begin{displaymath}
\mathbf{Q}=\begin{pmatrix}
    2 &1 &-3\\
    1 &5 &0\\
   -3 &0 &5
    \end{pmatrix}=\mathbf{L}\mathbf{L}^\top   \quad \mbox{with} \quad 
\mathbf{L}=\cfrac{1}{\sqrt{2}}\begin{pmatrix}
    2 &0 \\
    1 &3 \\ 
   -3 &1 \\
\end{pmatrix}.
\end{displaymath}
The Cholesky factor ${\bf L}$ above allows us to specifically write $p$ as the SOS
$$p({\bf x})=\cfrac{1}{2}(2x_1^2+x_1x_2-3x_2^2)^2+\cfrac{1}{2}(x_2^2+3x_1x_2)^2.$$
\smallskip

\section{Formulation of SOS polynomials}
The use of SOS programming in HJB problems have been addressed in prior works, e.g.\cite{lasagna2016sos}. In our work, we introduce three different SOS formulations:
\begin{itemize}  
     \item \textbf{Completing the sum of squares:} Adding higher degree terms to a given polynomial approximation (\ref{eq:pastEpoly}) or (\ref{eq:futureEpoly}) to make them sums of squares.
     \item \textbf{Completing the sum of squares with collocation:} Lower degree terms match a given polynomial approximation, and collocation is used to determine highest degree terms.
     \item \textbf{Collocation method:} Based on the Parrilo formulation of SOS (\ref{eq:psos}) within a least-squares collocation method.
\end{itemize}
The approaches are illustrated below using the past energy function $\mathcal{E}^-_\eta (\mathbf{x})$, but the same procedure can be applied to the future energy function  $\mathcal{E}^+_\eta (\mathbf{x})$ as well.  

\subsection{Completing the sum of squares}
Given a degree $d$ approximation to $\mathcal{E}_\eta^-(\mathbf{x})$ as in (\ref{eq:pastEpoly}), we propose a sum of squares approximation $\mathcal{E}_{sos}^-(\mathbf{x})$ as:
$$\mathcal{E}_{sos}^{-}(\mathbf{x})=\big( \mathbf{\tilde{v}}_1^\top  \mathbf{x}+\mathbf{\tilde{v}}_2^\top  \kronF{\mathbf{x}}{2}+\cdots +\mathbf{\tilde{v}}_{d-1}^\top  \kronF{\mathbf{x}}{d-1}\big)^2.$$
%
In general, the terms $\mathbf{\tilde{v}}_k$ are matrices and thus, this is a sum of squares.
The $d-1$ coefficients of $\mathcal{E}_{sos}^-(\mathbf{x})$ are found by matching the lowest degree terms in $\mathcal{E}_\eta^-(\mathbf{x})$. This is done without the involvement of the HJB equation. The HJB information is implicitly embedded in the polynomial approximation of $\mathcal{E}_\eta^-(\mathbf{x})$. \newline

\textbf{Example:} Given $n=2$, and a degree $d=3$ polynomial approximation $\mathcal{E}_\eta^-(\mathbf{x})$, the degree of the SOS approximation is $2(d-1)=4$. Therefore:
$$\mathcal{E}_{sos}^{-}(\mathbf{x})=\big( \mathbf{\tilde{v}}_1^\top  \mathbf{x}+\mathbf{\tilde{v}}_2^\top  \kronF{\mathbf{x}}{2}\big)^2,$$
where $\mathbf{\tilde{v}}_1^\top   \in \mathbb{R}^{2 \times 2}$ and $\mathbf{\tilde{v}}_2^\top   \in \mathbb{R}^{2 \times 4}$.  By expanding $\mathcal{E}^-_{sos}(\mathbf{x})$, we get:
$$\mathcal{E}_{sos}^{-}(\mathbf{x})= \mathbf{x}^\top  \mathbf{\tilde{v}}_1 \mathbf{\tilde{v}}_1^\top   \mathbf{x}+2\mathbf{x}^\top  \mathbf{ \tilde{v}}_1  \mathbf{\tilde{v}}_2^\top  \mathbf{\kronF{x}{2}}+(\mathbf{\kronF{x}{2}})^\top  \mathbf{\tilde{v}}_2 \mathbf{\tilde{v}}_2^\top   \mathbf{\kronF{x}{2}}$$ 
or
$$\mathcal{E}_{sos}^{-}(\mathbf{x})= {\tt vec}(\mathbf{\tilde{v}}_1 \mathbf{\tilde{v}}_1^\top  )^\top  \kronF{\mathbf{x}}{2}+2{\tt vec}( \mathbf{\tilde{v}}_1  \mathbf{\tilde{v}}_2^\top  )^\top  \kronF{\mathbf{x}}{3}+$$ $${\tt vec}(\mathbf{\tilde{v}}_2 \mathbf{\tilde{v}}_2^\top  )^\top \kronF{\mathbf{x}}{4}$$ 
Now, matching the $O(\kronF{\mathbf{x}}{2})$ and $O(\kronF{\mathbf{x}}{3})$ terms leads to the following systems of equations:
\begin{align*}
\mathbf{\tilde{v}}_1 \mathbf{\tilde{v}}_1^\top  &=\frac{1}{2} \mathbf{V}_2 \qquad \mbox{(solve by Cholesky)}\\
2\mathbf{\tilde{v}}_1 \mathbf{\tilde{v}}_2^\top  &=\frac{1}{2} \mathbf{V}_3 \qquad \mbox{(solve by backsubstitution)}
\end{align*}
where, $\mathbf{V}_2={\tt reshape}(\mathbf{v}_2^\top,2,2)$, $\mathbf{V}_3={\tt reshape}(\mathbf{v}_3^\top  ,2,4)$ and $\mathbf{v}_2^\top$ and $\mathbf{v}_3^\top$ are, respectively, the degree $2$ and $3$ coefficients of the polynomial approximation of $\mathcal{E}_\eta^-(\mathbf{x})$.  

\subsection{Completing the sum of squares with collocation}
In this approach, we write $\mathcal{E}_{sos,c}^-(\mathbf{x})$ in the same polynomial squared form with an additional term:
$$\mathcal{E}_{sos,c}^{-}(\mathbf{x})\equiv\big( \mathbf{\tilde{v}}_1^\top  \mathbf{x}+\mathbf{\tilde{v}}_2^\top  \kronF{\mathbf{x}}{2}+\cdots +\mathbf{\tilde{v}}_{d-1}^\top  \kronF{\mathbf{x}}{d-1}+{\mathbf{\tilde{v}}_d^\top  }\kronF{\mathbf{x}}{d}\big)^2$$
and we calculate for the first $d-1$ coefficients as above (matching to a degree $2d-1$ polynomial approximation).

To calculate $\mathbf{\tilde{v}}_d$ we introduce the residual of the HJB equation found by substituting $\mathcal{E}_{sos,c}^-(\mathbf{x})$ into (\ref{eq:HJB-NLHinfty1}) at the point $\mathbf{x}$ which we denote by $\mathbf{R}(\mathbf{x})$.  The least-squares collocation problem is then described by the functional
$$J({\mathbf{\tilde{v}}_d^\top  })=\sum_{k=1}^N \|\mathbf{R}(\mathbf{x}^{(k)})\|^2,$$
that is formed using $N$ sample points $\{\mathbf{x}^{(k)}\}$ in the approximation domain $\Omega \subset \mathbb{R}^n$.  By using randomly selected sample points within $\Omega$, this procedure has the potential to scale to larger dimensions. The missing coefficient in $\mathcal{E}_{sos,c}^-(\mathbf{x})$ is then given by ($*$ in this work denotes optimal not an adjoint)
$${\tilde{\mathbf{v}}_d^{\top *}  }=\,\mathop{\mbox{argmin}}_{\tilde{\mathbf{v}}_d^\top   \in \mathbb{R}^{n^d}} \,J(\tilde{\mathbf{v}}_d^\top  ). $$
 
\subsection{Nonlinear least-squares collocation\label{sec:SOS}}
 
Using Parrilo's formulation in (\ref{eq:psos}), we write $\mathcal{E}_{p}^{-}(\mathbf{x})$ as:
\begin{equation}
\label{eq:sos}
\mathcal{E}_{p}^{-}(\mathbf{x})=\mathbf{z}(\mathbf{x})^\top  \mathbf{Q}\mathbf{z}(\mathbf{x})   
\end{equation}
where $\mathbf{z}(\mathbf{x}) \in \mathbb{R}^{\nu}$ is a vector of all monomials of degree $\leq d$ and $\mathbf{Q} \in \mathbb{R}^{\nu \times \nu}$ is a positive definite matrix. The number of monomial terms, $\nu$, in $\mathbf{z}(\mathbf{x})$ is defined as follows:
$$\nu= \sum_{i=1}^d {\rm deg}_i(n),$$
where ${\rm deg}_i(n)$ is the number of monomial terms of degree $i$ in dimension $n$ defined recursively by:
$${\rm deg}_i(n)=\sum_{j=1}^n {\rm deg}_{i-1}(j) \quad \text{for}\,\, i\geq 2$$
with ${\rm deg}_1(n)=n$.  The matrix
$\mathbf{Q}$ can be factorized as $\mathbf{Q}=\mathbf{L}\mathbf{L}^\top  $ and $\mathcal{E}_p^-(\mathbf{x})$ can be rewritten as: 
\begin{equation}\label{eq:psosL}
  \mathcal{E}_{p}^{-}(\mathbf{x})=p_{sos}(\mathbf{x})=\mathbf{z}(\mathbf{x})^\top  \mathbf{L}\mathbf{L}^\top  \mathbf{z}(\mathbf{x})
\end{equation}
where the entries of $\mathbf{L}$ are the decision variables of the optimization problem defined by the objective function 
\begin{equation}\label{eq:Objfun}
 J(\mathbf{L})=\sum_{k=1}^N \|{\bf R}(\mathbf{x}^{(k)})\|^2,   
\end{equation}
where now ${\bf R}(\mathbf{x})$ is the residual of the HJB equation at $\mathbf{x}$ when substituting (\ref{eq:psosL}) into (\ref{eq:HJB-NLHinfty1}).  The optimal coefficients in (\ref{eq:psosL}) are then given by
$${\mathbf{L}^*}=\,\mathop{\mbox{argmin}}_{\mathbf{L} \in \mathbb{R}^{\nu \times \nu }} \,J(\mathbf{L}). $$

\subsection{Sequence of polynomial approximations}
In these last two sections, we performed an optimization problem to compute unknown coefficients over the domain $\Omega$.  In the nonlinear least-squares collocation case especially, the entries of $\mathbf{L}$ strongly depend on the region of approximation. For approximation regions $\Omega \subset \mathbb{R}^n$ that lie within the unit hypercube, the lowest degree terms in the residual ${\bf R}$ tend have a more significant contribution than the higher degree terms. On the other hand, for approximation regions $\Omega$ where most sample points lie outside the unit hypercube, the higher degree terms have more contribution the residual $J(\mathbf{L})$. This motivates a sequence of polynomial approximations that increase as the approximation domain $\Omega$ increases, and has the effect of warm-starting the optimization algorithm.\\

\textbf{Example:} Consider the case $d=4$ where the system has a cubic nonlinearity $\dot{\mathbf{x}}=\mathbf{A}\mathbf{x}+\mathbf{F}\kronF{\mathbf{x}}{3}+\mathbf{B}\mathbf{u}$. The energy function takes the following form 
$$\mathcal{E}_{p}(\mathbf{x})=\big[\mathbf{x}^\top  {\kronF{\mathbf{x}}{2}}^\top   \big]
\begin{pmatrix}
\mathbf{L}_{11} & 0 \\
\mathbf{L}_{21} & \mathbf{L}_{22} 
\end{pmatrix}
\begin{pmatrix}
\mathbf{L}_{11} & 0 \\
\mathbf{L}_{21} & \mathbf{L}_{22} 
\end{pmatrix}^\top  
\big[\mathbf{x}^\top {\kronF{\mathbf{x}}{2}}^\top   \big]^\top.  $$
We can expand this to find
\begin{IEEEeqnarray*}{lCr}
\mathcal{E}_{p}(\mathbf{x})&=&\mathbf{x}^\top  \mathbf{L}_{11}\mathbf{L}_{11}^\top\mathbf{x}+2\mathbf{x}^\top  \mathbf{L}_{11}\mathbf{L}_{21}^\top{\kronF{\mathbf{x}}{2}}\\
&& +{\kronF{\mathbf{x}}{2}}^\top  (\mathbf{L}_{12}\mathbf{L}_{12}^\top
+\mathbf{L}_{22}\mathbf{L}_{22}^\top){\kronF{\mathbf{x}}{2}}.    
\end{IEEEeqnarray*}
We let $\mathbf{Q}_1=\mathbf{L}_{11}\mathbf{L}_{11}^\top  $, $\mathbf{Q}_2=\mathbf{L}_{11}\mathbf{L}_{21}^\top  $ and $\mathbf{Q}_3=\mathbf{L}_{21}\mathbf{L}_{21}^\top+\mathbf{L}_{22}\mathbf{L}_{22}^\top  $. Note that $\mathbf{Q}_3$ is the only matrix that depends on $\mathbf{L}_{22}$.
Substituting $\mathcal{E}_{p}(\mathbf{x})$ into (\ref{eq:HJB-NLHinfty1}), the residual ${\bf R}(\mathbf{x})$ is a degree 6 polynomial,  $\mathbf{Q}_3$, and hence $\mathbf{L}_{22}$, only appears in the degree 4, 5, and 6 terms, $r_4(\mathbf{L})$, $r_5(\mathbf{L})$ and $r_6(\mathbf{L})$, respectively.  These coefficient functions are
\begin{align*}
 r_4(\mathbf{L})=&{\tt vec}(8\mathbf{Q}_1\mathbf{B}\mathbf{B}^\top  \Tilde {\mathbf{Q}}_3^\top  +2\mathbf{Q}_1\mathbf{F}+4\mathbf{A}^\top  \Tilde {\mathbf{Q}}_3^\top  )\\
 &{ }+ {\tt vec}(18\Tilde {\mathbf{Q}}_2 \mathbf{B}\mathbf{B}^\top  \Tilde {\mathbf{Q}}_2^\top  )\\
r_5(\mathbf{L})=&{\tt vec}(24\Tilde {\mathbf{Q}}_2\mathbf{B}\mathbf{B}^\top  \Tilde {\mathbf{Q}}_3^\top  +6\Tilde {\mathbf{Q}}_2\mathbf{F})\\
r_6(\mathbf{L})=&{\tt vec}(8\Tilde {\mathbf{Q}}_3\mathbf{B}\mathbf{B}^\top  \Tilde {\mathbf{Q}}_3^\top  +4 \Tilde {\mathbf{Q}}_3\mathbf{F}), 
\end{align*}
where $\Tilde {\mathbf{Q}}_2$ and $\Tilde {\mathbf{Q}}_3$ are reshaped versions of 
$\mathbf{Q}_2$ and $\mathbf{Q}_3$ to match the dimensions.   When we attempt to minimize the residual with many sample points far outside the unit hypercube ($\| \mathbf{x}\|\gg 1$), the degree 6 term dominates and therefore the objective function is very sensitive to $\mathbf{L}_{22}$.  Small changes to ${\bf \tilde{Q}}_3$ have such a large influence that the $r_6({\bf L})$ term dominates the residual and this tends to drive ${\bf \tilde{Q}}_3$ to zero, but this also impacts ${\bf L}_{21}$ using the definition of ${\bf Q}_3$ above.  Therefore, we elect to drop the $\mathbf{L}_{22}$ term in the case $\Omega$ is far outside the unit hypercube. 
In the case of more blocks in ${\bf L}$, it is desirable to drop the term representing the highest degree coefficients. For example, for degree 6 SOS approximations, we only keep the shaded portion of the matrix ${\bf L}$ illustrated below when $\Omega$ has more sample points outside the unit hypercube:
\begin{center}   
\begin{tikzpicture}
  \matrix[matrix of math nodes,left delimiter = (,right delimiter = ),row sep=7pt,column sep = 7pt] (m)
  {
    L_{11} & 0      &  0\\
    L_{21} & L_{22} &  0 \\
    L_{31} & L_{32} & L_{33} \\
  };
  \begin{pgfonlayer}{background}
    \node[inner sep=3pt,fit=(m-1-1)]          (1)   {};
    \node[inner sep=3pt,fit=(m-2-1) (m-3-2)]  (2)   {};
    \node[inner sep=3pt,fit=(m-3-3)]          (3)   {};
    \draw[rounded corners,dotted,fill=green!50!white,inner sep=3pt,fill opacity=0.1]
    (1.north west) |- (2.south east) |- (2.east) |- (2.north) |- (1.north) -- cycle;
  \end{pgfonlayer}
\end{tikzpicture}
\end{center}

\subsection{Windowing and issues for large regions $\Omega$ in $\mathbb{R}^n$ }
\begin{figure}
    \centering
    \includegraphics[scale=0.4]{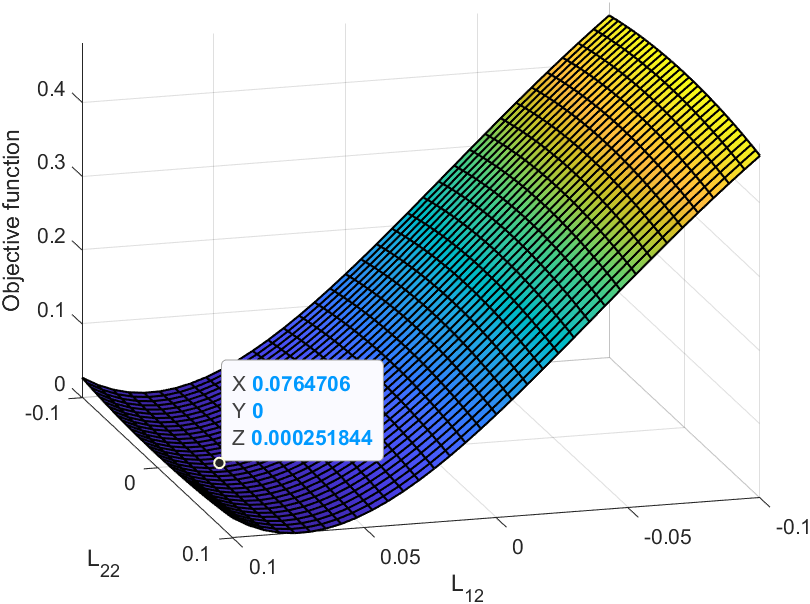}  
   \caption{Objective function when fixing $\mathbf{L}_{11}$ at $\mathbf{L}_{11}^*$ and choosing collocation points in the range $[-1, 1]$.}
   \label{fig1}
\end{figure}

\begin{figure}
    \centering
    \includegraphics[scale=0.4]{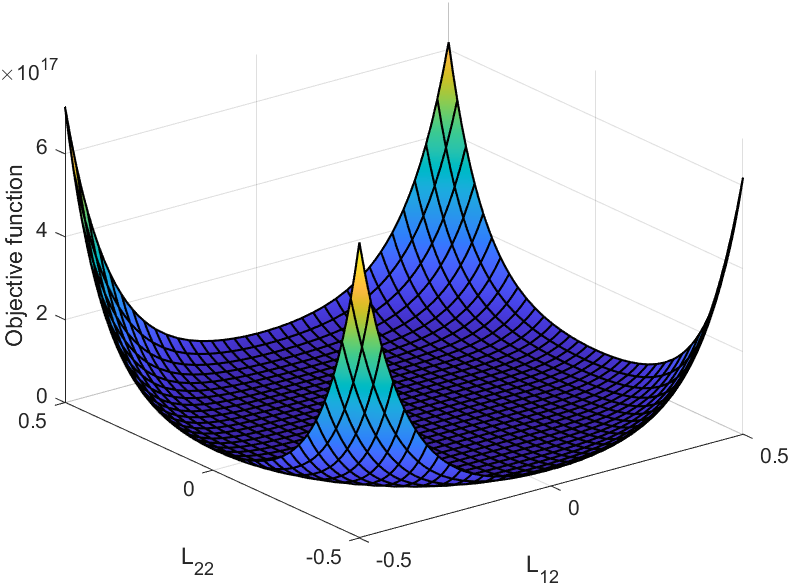}  
   \caption{Objective function when fixing $\mathbf{L}_{11}$ at $\mathbf{L}_{11}^*$ and choosing collocation points in the range $[-20, 20]$.}
   \label{fig2}
\end{figure}

\begin{figure}
    \centering
    \includegraphics[scale=0.4]{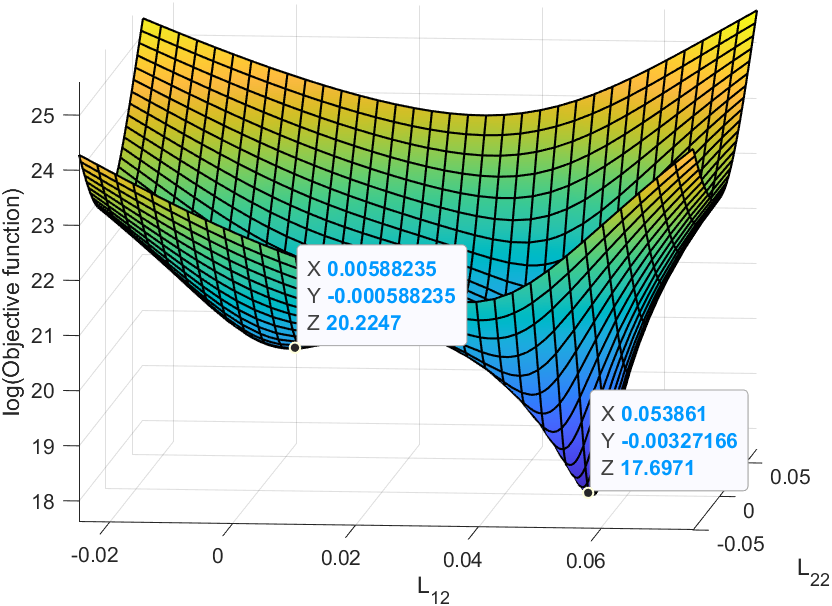}  
   \caption{Log of the objective function when fixing $\mathbf{L}_{11}$ at $\mathbf{L}_{11}^*$ and choosing collocation points in the range $[-20, 20]$.}
   \label{fig3}
\end{figure}

Performing the optimization of $J({\bf L})$ for small regions around the origin is typically a convex optimization problem. To illustrate this, we take the scalar system from \cite{kramer2023nonlinear1} and we write the degree $4$ SOS approximation as in (\ref{eq:psosL}) for the past energy function ${\mathcal{E}}_p^{-}({\mathbf{x}})$ with  an approximation region $\Omega=[-1;1] \subset \mathbb{R}$. Let $\mathbf{L}_{11}^*$ be the Cholesky factor of  the solution of the Riccati equation (\ref{eq:12}). Now, if we fix the $\mathbf{L}_{11}$ entry of $\mathbf{L}$ at $\mathbf{L}_{11}^*$, Fig.~\ref{fig1} represents the objective function with respect to $\mathbf{L}_{12}$ and $\mathbf{L}_{22}$ (the {\em insensitivity} of $J$ to $\mathbf{L}_{22}$ when most samples are {\em inside} the unit interval is clearly illustrated in this figure as well).
If we consider larger domains, e.g. $\Omega=[-20;20] \subset \mathbb{R}$, Figs.~\ref{fig2} and \ref{fig3} show that the optimization problem quickly becomes non-convex.

This non-convexity can be partially overcome by introducing a {\em windowing} procedure that warm starts a sequence of optimization problems over progressively larger regions $\Omega_1\subset\cdots\subset\Omega_k\subset \cdots\Omega$.
We start a collocation problem for $\Omega_1$, a small region around the origin, where the SOS approximation will easily capture the behaviour of the lower order terms of the polynomial approximation. Then, the optimization results for ${\bf L}$ are used as warm starts (initial guesses) for problems with successively larger domains $\Omega_k$, while we also increase the number of sample points to $s_k$.   In this way, we achieve     
$$\mathbf{L}^{*(s_1)}_{\Omega_1} \rightarrow \mathbf{L}^{*(s_2)}_{\Omega_2} \rightarrow \cdots \rightarrow \mathbf{L}^{*(N)}_{\Omega}.$$
\subsection{SOS nonlinear least-squares collocation algorithm\label{sec:collocation}}
\textbf{Algorithm 1} summarizes the procedure of approximating energy functions using the SOS collocation method described in Section~\ref{sec:SOS}.
We begin by solving the algebraic Riccati equation corresponding to the type of energy function we want to approximate, (\ref{eq:12}) or (\ref{eq:13}) to provide an initial guess for the ${\bf L}_{11}$ term. Then, we set the appropriate dimension $\nu$ of the $\mathbf{L}$ matrix: depending on the dimension of the system $n$, the degree of the approximation $d$, and the sequence of polynomial approximations (deciding on diagonal blocks of $\mathbf{L}$ to keep. The optimization problem is solved in a local region around the origin, to get the first warm start on $\Omega_1$.  We then proceed with windowing until the desired approximation region is reached. Finally, we construct the $\mathbf{Q}=\mathbf{L}\mathbf{L}^\top$ matrix that defines the energy function.\\

\noindent\fbox{%
        \parbox{.95\columnwidth}{%
        \textbf{Algorithm 1:} SOS approximation of energy functions. \\\\
        \textbf{Input:} System matrices $\mathbf{A}$, $\mathbf{B}$, $\mathbf{C}$, $\mathbf{F}$, $\eta$, $d$ the degree of SOS approximation, a nested set of approximation regions $\Omega_0 \subset \Omega_1 \cdots \subset \Omega_r= \Omega$  and corresponding numbers of sample points $\{s_i\}_{i=0}^r$ ($s_r=N$).\\
        \textbf{Output:} $\mathbf{Q}$, where $\mathbf{Q}=\mathbf{L}\mathbf{L}^\top$ in (\ref{eq:sos}).
        \begin{enumerate}

       \item Solve the Algebraic Riccati equation (\ref{eq:12}) or (\ref{eq:13}), corresponding to the desired energy function, to get a warm start for the lowest degree term.
       \item Calculate $\nu$ that depends on both $n$ and $d$.
       \item Calculate, $\nu_1$, the appropriate dimension that depends on the structure of the $\mathbf{L}$ resulting from the sequence of polynomial approximations.
       \item To get a warm start $\mathbf{L}_0$ for the optimization, solve:
       $$\mathbf{L}_{0}=\,\mathop{\mbox{argmin}}_{\mathbf{L} \in \mathbb{R}^{\nu \times \nu_1 }} \,\sum_{k=1}^{s_0} \|{\bf R}(\mathbf{x}^{(k)})\|^2$$
       where $\{{\bf x}^{(k)}\}_{k=1}^{s_0}\subset\Omega_0$, a small region containing ${\bf 0}$.
       \item 
       For $i=1,2,\cdots,r$
       
          \ Solve: $$\mathbf{L}_{i}=\,\mathop{\mbox{argmin}}_{\mathbf{L} \in \mathbb{R}^{\nu \times \nu_1 }} \,\sum_{k=1}^{s_i} \|{\bf R}(\mathbf{x}^{(k)})\|^2$$
          \ with $\mathbf{L}_{i-1}$ as a warm start, and
          
          \ $s_i$ is the number of samples  in $\Omega_i$. 
       
       end \\
       
           \item Finally, set $\mathbf{L}=\mathbf{L}_r$ and compute $\mathbf{Q}=\mathbf{L}\mathbf{L}^\top$.
              
       \end{enumerate} 
       }
 }%

\section{Numerical Results}
We will perform numerical demonstrations of the proposed methods using a scalar example\cite{kramer2023nonlinear1} to approximate the past energy $\mathcal{E}_\eta^-(\mathbf{x})$.  Then we will test Algorithm 1 on a 6-dimensional ring of van der Pol oscillators where we approximate the future energy function $\mathcal{E}_\eta^+(\mathbf{x})$.  Finally, we consider a 12-dimensional problem resulting from finite element approximation to Burgers equation. For the optimization, we used the IPOPT algorithm~\cite{IPOPT} along with the Julia package JuMP~\cite{Lubin2023}. 

\subsection{Scalar example}
\begin{figure}[!ht]
   \centering
   \includegraphics[scale=0.4]{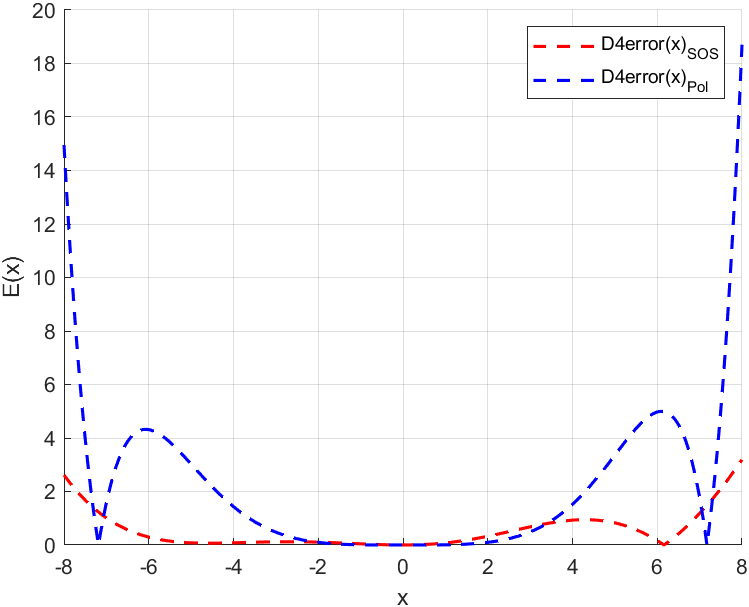}  
   \caption{Error of degree 4 approximations of energy functions $[-8,8]$}
   \label{fig4}
   \end{figure}

\begin{figure}[!ht]  
\centering
   \includegraphics[scale=0.4]{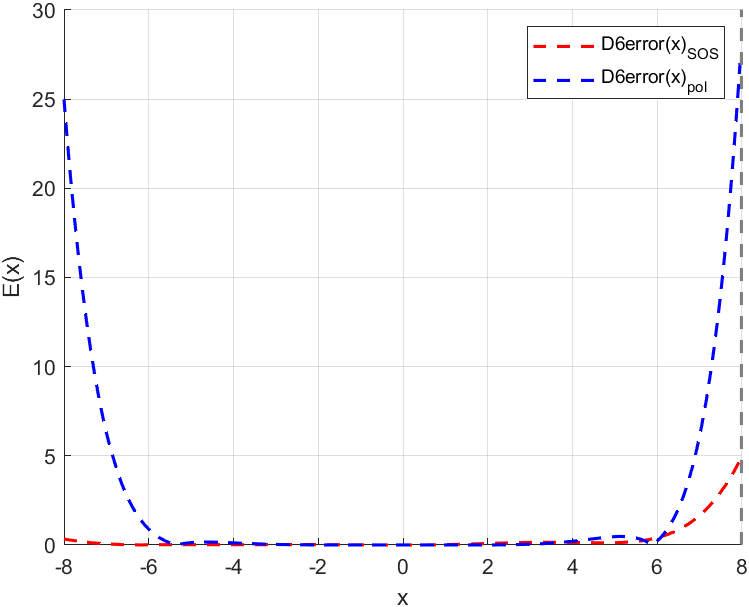}  
   \caption{Error of degree 6 approximations of energy functions $[-8,8]$}
   \label{fig5}
\end{figure}

\begin{figure}[!ht]
\centering
   \includegraphics[scale=0.4]{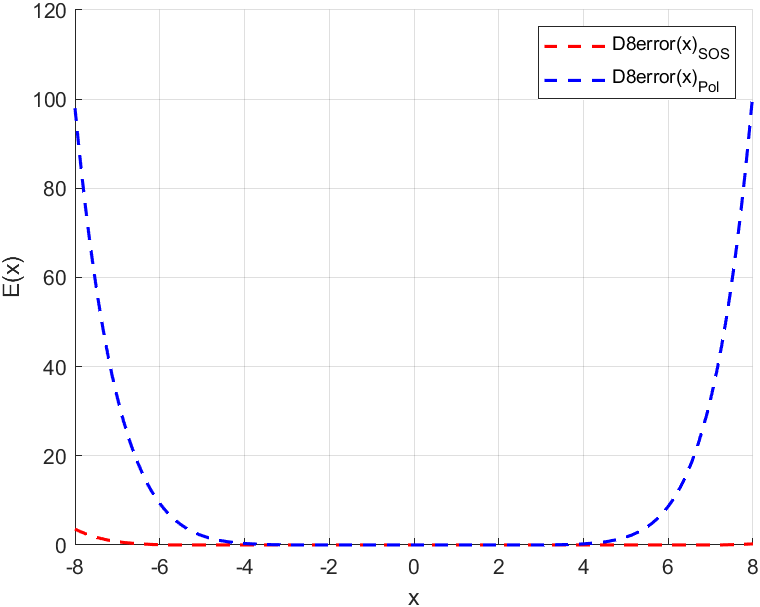}  
  \caption{Error of degree 8 approximations of energy functions $[-8,8]$}
  \label{fig6}
\end{figure}

For our first study, we approximate the past energy function $\mathcal{E}_{\eta}^-(\mathbf{x})$ ($\gamma=\sqrt{2}$, $\eta=0.5$) with degree 4, 6, and 8 SOS polynomials, for a scalar example found in \cite{kramer2023nonlinear1} using $\Omega=[-8,8]$. The system is given by
\begin{displaymath}
  \dot{x}(t)=-2x(t)+x^2(t)+2u(t) \quad \mbox{and} \quad y(t)=2x(t).
\end{displaymath}
We used symmetric windows starting at $[-1,1]$ and doubling the size until we reached $[-8,8]$. Since we can solve the HJB in the 1D setting, the results of the polynomial approximations in Section~\ref{sec:polynomial} and \ref{sec:collocation} are compared to the analytical solution. The errors in the polynomial approximation, and the SOS approximation are presented in Figs.~\ref{fig4}, \ref{fig5}, and \ref{fig6}. 

As we increase the degree of the SOS approximation, the overall error gets smaller and is superior to the polynomial approximation in \cite{kramer2023nonlinear1}. As expected, the polynomial approximation is good locally, but has a poor global behaviour that does not preserve the positivity and monotonicity properties that are expected for this example. 

\subsection{Ring of van der Pol equations}
Consider the control of a ring of van der Pol oscillators,
\begin{equation}\label{eq:vdp}
  \ddot{y}_i + (y_i^2-1)\dot{y}_i + y_i = y_{i-1}-2y_i+y_{i+1} + b_i u_i(t),
\end{equation} 
for $i=1,\ldots,g$ with $y_i(0)=y^0$ and $\dot{y}_i(0)=0$ (we identify $y_{g+1}=y_1$ and $y_g=y_0$ to 
close the ring).  The coupling terms on the right-hand-side resemble a discrete Laplacian. 
The stability of this system
was studied in \cite{nana2006SynchronizationRingFour} and a related control problem considered in \cite{barron2016StabilityRingCoupled}.  Here, we will consider the case when $g=3$, with $b_1=b_2=1$, $b_3=0$.  Writing this as a first-order system places the system in the form (\ref{eq:nl_system}) with $n=6$ states and $m=2$ control inputs. We define $\mathbf{C}$ using the positions of the 3 oscillators. We demonstrate our proposed SOS approximation on the future energy function $\mathcal{E}^+_{\eta}(\mathbf{x})$ using a degree 4 SOS approximation for the case $\eta=1$. 
   
Since we do not have an analytical solution available, we will evaluate the quality of the approximation by deriving the optimal control for the system using both the SOS approximation and the polynomial approximations
\begin{equation}
\label{eq:OptCont}
\begin{split}
    &u^*_{sos}(t)=-\mathbf{B}^\top  \cfrac{\partial \mathcal{E}_{sos}^{+}(\mathbf{x})}{\partial x}\\  
    &u^*_{pol}(t)=-\mathbf{B}^\top  \cfrac{\partial \mathcal{E}_{pol}^{+}(\mathbf{x})}{\partial \mathbf{x}}
\end{split}
\end{equation}
running the closed-loop simulations of (\ref{eq:vdp}) and simultaneously approximating the integral term in the energy function given as  
\begin{equation}
\label{eq:Integral}
 {\mathcal{E}}_\eta^{+}(\mathbf{x}_0)= \min_{\substack{u\in L_2[0,\infty)\\ \mathbf{x}(-\infty)=0\ {\mathbf{x}}(0)=\mathbf{x}_0}} \cfrac{1}{2}\int^{\infty}_0||y(t)||^2+||u(t)||^2 dt,   
\end{equation}
where $x_0$ is the initial condition for the system.
We used 5 symmetric windows starting at $[-0.1, 0.1]^6$ and increasing the size by 0.1 in each direction until we reached $[-0.5, 0.5]^6$. To remove bias in our tests, we performed study using $1000$ random initial conditions $\mathbf{x}_0$ in the approximation domains. The number of times the closed-loop system is unstable (or if was always stable) is reported in Table~\ref{tab:vdp_stability}.
\begin{table} 
\centering
\caption{\label{tab:vdp_stability}Stability of closed-loop systems from sampled initial guesses.}
\begin{tabular}{ |c|c|c| } 
 \hline
 Windows in $\mathbb{R}^{6}$ & Poly Approx & SOS Approx \\ 
 \hline\hline
 $[-0.1,\,0.1]^6$& Stable & Stable \\  
 \hline
 $[-0.2,\,0.2]^6$& Stable & Stable  \\
 \hline
 $[-0.3,\,0.3]^6$& Unstable 5/1000 & Stable  \\
 \hline
 $[-0.4,\,0.4]^6$& Unstable 51/1000& Stable  \\
 \hline
 $[-0.5,\,0.5]^6$& Unstable 116/1000 & Stable  \\ 
 \hline
\end{tabular}
\end{table}
To show that the sum of squares polynomial approximation has comparable accuracy, we report the average relative errors in Table~\ref{tab:vdp_relErr}.  Note that we did not include the unstable results from the polynomial approximations in the reported averages (since the simulations didn't finish).
\begin{table} 
\centering
\caption{\label{tab:vdp_relErr}Average relative errors in $\mathcal{E}_1^+$ from 1000 initial conditions.  The errors from the unstable cases were not included.}
\begin{tabular}{ |c|c|c| } 
 \hline
 Windows in $\mathbb{R}^{6}$ & Poly Approx error & SOS Approx error \\ 
 \hline\hline
 $[-0.1,\,0.1]^6$ & $3.5465 \times 10^{-4}$ & $4.1133 \times 10^{-3}$ \\  
 \hline
 $[-0.2,\,0.2]^6$ & $5.1609 \times 10^{-4} $ & $1.5224 \times 10^{-2}$  \\
 \hline
 $[-0.3,\,0.3]^6$ & $2.9679 \times 10^{-3}$ & $3.0859\times 10^{-2}$  \\
 \hline
 $[-0.4,\,0.4]^6$ & $ 1.2120\times 10^{-2}$ & $5.1790 \times 10^{-2}$  \\
 \hline
 $[-0.5,\,0.5]^6$ & $ 1.7436\times 10^{-2}$ & $7.1467 \times 10^{-2}$  \\
 \hline
\end{tabular}
\end{table}

 Table~\ref{tab:vdp_stability} shows that the SOS approximation provided a stable closed-loop solution for each of the 1000 random starting points in each approximation domain. On the other hand, as we increased the size of approximation region the polynomial approximation from \cite{kramer2023nonlinear1} tends to give unstable solutions.  This is to be expected since the polynomial approximations are designed to be very accurate near the origin, but without being an SOS polynomial, there is no guarantee that positivity is maintained away from the origin.  For example, starting at the initial condition,
 \begin{displaymath}
   \mathbf{x}_0^{\text{unstable}}=\begin{pmatrix}
 -0.21 &  0.08 &  0.06& -0.35&  0.36& -0.47
 \end{pmatrix}^\top
\end{displaymath}
leads to an unstable simulation. Table~\ref{tab:vdp_relErr} shows that the average relative error, excluding the unstable initial conditions for the polynomial approximation, is acceptable for both the SOS and the polynomial approximation, with a clear advantage for the polynomial approximation in each of the approximation regions.  The main advantage of the SOS approximations is that every initial condition in the neighborhood of the origin can be stabilized using (\ref{eq:OptCont}).
\subsection{Burgers equation}
We also consider the control of the discretized 1D Burgers equation (using linear finite elements).  Consider
\begin{equation}
\label{eq:Burgers}
z_t = -z z_x + \epsilon z_{xx}= \sum_{i=1}^m \chi_i^m(x)u_i(t), \,\, x\in(0,1),\,\, t>0, 
\end{equation} 
with periodic boundary conditions and controlled outputs $y_i(t) = \int_0^1 \chi_i^p(x)z(t,x)\ dx$, $i=1,\ldots,p$.  Here, $\chi_i^r(x)$ is the characteristic function over $[(i-1)/r,i/r]$ and $\epsilon=5\times 10^{-3}$.  Approximating this PDE using 12 linear finite elements and choosing $m=p=6$ creates a quadratic system of equations ($n=12$) of the form 
\begin{equation} 
\label{eq:DisBurgEqu}
\dot{\mathbf{x}}=\mathbf{Ax+N}\kronF{\mathbf{x}}{2}+\mathbf{B}u, \qquad \mbox{and} \qquad y=\mathbf{Cx}.
\end{equation} 
We demonstrate our proposed SOS approximation on the future energy function $\mathcal{E}_\eta^+(\mathbf{x})$
using a degree $4$ SOS approximation for the case $\eta=1$. As for the ring of van der Pol equations,  we do not have an analytical solution. We will evaluate the quality of the approximation by deriving the optimal control for the discritized Burgers equation using both the SOS approximation and the polynomial approximation (\ref{eq:OptCont}) and running the closed-loop simulations of (\ref{eq:DisBurgEqu}) and simultaneously approximating the integral term in the energy function given by (\ref{eq:Integral}). For this example we ran the optimization for one window $[-0.1,0.1]^{12}$, with 1000 random initial conditions $\mathbf{x}_0$ in the same window. To show the effectiveness of the approximation outside the optimization range, we also take 1000 random initial conditions $\mathbf{x}_0$ in $[-0.2,0.2]^{12}$, $[-0.3,0.3]^{12}$ and $[-0.4,0,4]^{12}$ outside of the approximation window . We report the average relative errors in Table \ref{tab:Burgers_relErr}. 
The results in Table \ref{tab:Burgers_relErr} show that in this case both the SOS approximation and the polynomial approximation provide a good optimal control for the different random starting values both in the SOS approximation window and outside of it. 
\begin{table} 
\centering
\caption{\label{tab:Burgers_relErr}Average relative errors in $\mathcal{E}_1^+$ from 1000 initial conditions.  The errors from the unstable cases were not included.}
\begin{tabular}{ |c|c|c| } 
 \hline
 Windows in $\mathbb{R}^{12}$ & Poly. approx. error & SOS approx. error \\ 
 \hline\hline
 $[-0.1,\,0.1]^{12}$ & $7.3851 \times 10^{-4}$ & $3.4814\times 10^{-2}$  \\  
 \hline
 $[-0.2,\,0.2]^{12}$ & $8.1425 \times 10^{-4}$ & $4.5553\times 10^{-2}$  \\  
 \hline
  $[-0.3,\,0.3]^{12}$ & $1.0639 \times 10^{-3}$ & $6.0839\times 10^{-2}$  \\  
 \hline
  $[-0.4,\,0.4]^{12}$ & $1.8511 \times 10^{-3}$ & $7.5802\times 10^{-2}$  \\  
 \hline
\end{tabular}
\end{table}

\section{Conclusions}
Three sum of squares polynomial approximation methods were described.  In practice, the approach based on the form provided by \cite{parillo2000structured} lead to the best results.  The superiority of the SOS approximation over a Taylor series-based approximation was demonstrated using a scalar example with a known analytic formula.  In a more challenging case, the results were mixed.  The SOS approximation always provided stable closed-loop systems, but the Taylor polynomial approximations were more accurate on average when they produced stable closed-loop systems.  However, away from the origin, the Taylor polynomial approximations failed to generate a stabilizing control in more than 10\% of the randomly generated cases while the SOS polynomial approximations produced stable closed-loop systems in each of those same cases.  The challenge of solving the nonlinear collocation problem limits the scalability of this approach to more modest sizes and higher degree SOS polynomials due to the large number of monomial terms that appear in the $\mathbf{L}$ matrices.  Further work is underway to find better formulations of the optimization problems and test the performance on larger systems.

\section*{Acknowlgement}
This work was supported in part by the National Science Foundation under Grant CMMI-2130727.

\bibliographystyle{plain}
\bibliography{references}

\end{document}